\documentclass[12pt]{article}
\usepackage{amssymb}
\usepackage{amsmath,amsthm}
\usepackage[latin1]{inputenc}
\usepackage{graphicx}
\usepackage{epsfig}
\usepackage{hyperref}
\hypersetup{colorlinks=true, urlcolor=blue}

 \newtheorem{remark}{Remark}

 \newtheorem{theorem}[remark]{Theorem}

\title{Functional centrality in graphs}

\author{\href{http://deim.urv.cat/~jarodriguez/}{J. A. Rodr\'{\i}guez}\footnote{e-mail:\mbox{\tt
    juanalberto.rodriguez\@@urv.cat}} \\
{\em Department of Computer Engineering and Mathematics}\\
Rovira i Virgili University of Tarragona\\ Av. Pa\"{\i}sos Catalans
26, 43007 Tarragona, Spain\\
E. Estrada\footnote{estrada66@yahoo.es}
\\ {\em Complex Systems Research Group} \\
 X-rays Unit, RIAIDT, Edificio CACTUS \\
  University of Santiago de Compostela,
\\15706 Santiago de Compostela, Spain\\
A. Gutiérrez\footnote{e-mail:\mbox{\tt
    amauri@bayesinf.com}} \\
{\em Bayes Inference, S. A.}\\ Gran Vía 39, 5º planta \\
28013 Madrid, Spain}

\date{}

\begin{document}

\maketitle

\begin{abstract}
In this paper we introduce the functional centrality as a
generalization of the subgraph centrality. We propose a general
method for characterizing nodes in the graph according to the number
of closed walks starting and ending at the node. Closed walks are
appropriately weighted according to the topological features that we
need to measure.
\end{abstract}



\section{Introduction}

A kind of local characterization of networks is made numerically by
using one of several measures known as \emph{centrality}
\cite{freeman}. One of the most used centrality measures is the
\emph{degree centrality}, DC \cite{Barabasi2}, which is a
fundamental quantity describing the topology of scale-free networks
\cite{SmallW}. DC can be interpreted as a measure of immediate
influence, as opposed to long-term effect in the network
\cite{freeman}. For instance, if a certain proportion of nodes in
the network are infected, those nodes having a direct connection
with them will also be infected. However, although a node in a
network may be linked to only one node, the risk of infection to the
first node remains high if the latter is connected to many others.

There are several other centrality measures that have been
introduced and studied for real world networks, in particular for
social networks. They account for the different node characteristics
that permit them to be ranked in order of importance in the network.
\emph{Betweenness centrality} (BC) characterizes how influential a
node is in communicating between node pairs \cite{Freeman2}. In
other words, BC measures the number of times that a shortest path
between nodes $i$ and $j$ travels through a node $k$ whose
centrality is being measured. The farness of a vertex is the sum of
the lengths of the geodesics to every other vertex. The reciprocal
of farness is \emph{closeness centrality} (CC). The normalized
closeness centrality of a vertex is the reciprocal of farness
divided by the minimum possible farness expressed as a percentage
\cite{Barabasi2,freeman}. This measure is only applicable to
connected networks, since the distance between unconnected nodes is
undefined. Neither BC nor CC can be related to the network subgraphs
in a way that permits them to be considered as measures of node
subgraph centrality.

A centrality measure that is not restricted to shortest paths is the
\emph{eigenvector centrality} (EC) \cite{Bonacich}, which is defined
as the principal or dominant eigenvector of the adjacency matrix A
representing the connected subgraph or component of the network. It
simulates a mechanism in which each node affects all of its
neighbors simultaneously \cite{Bonacich2}. EC cannot be considered
as a measure of centrality whereby nodes are ranked according to
their participation in different network subgraphs. For instance, in
a graph with all nodes having the same degree (a regular graph), all
the components of the main eigenvalue are identical
\cite{Cvetkovic}, even if they participate in different subgraphs.
EC is better interpreted as a sort of extended degree centrality
which is proportional to the sum of the centralities of the node'
neighbors. Consequently, a node has high value of EC either if it is
connected to many other nodes or if it is connected to others that
themselves have high EC \cite{Newman}.


\begin{figure}[h]\label{fig}
\begin{center}
\vspace{4cm}
\begin{picture}(1,1)
\hspace{-6cm}\includegraphics[angle=90, width=5.5cm]{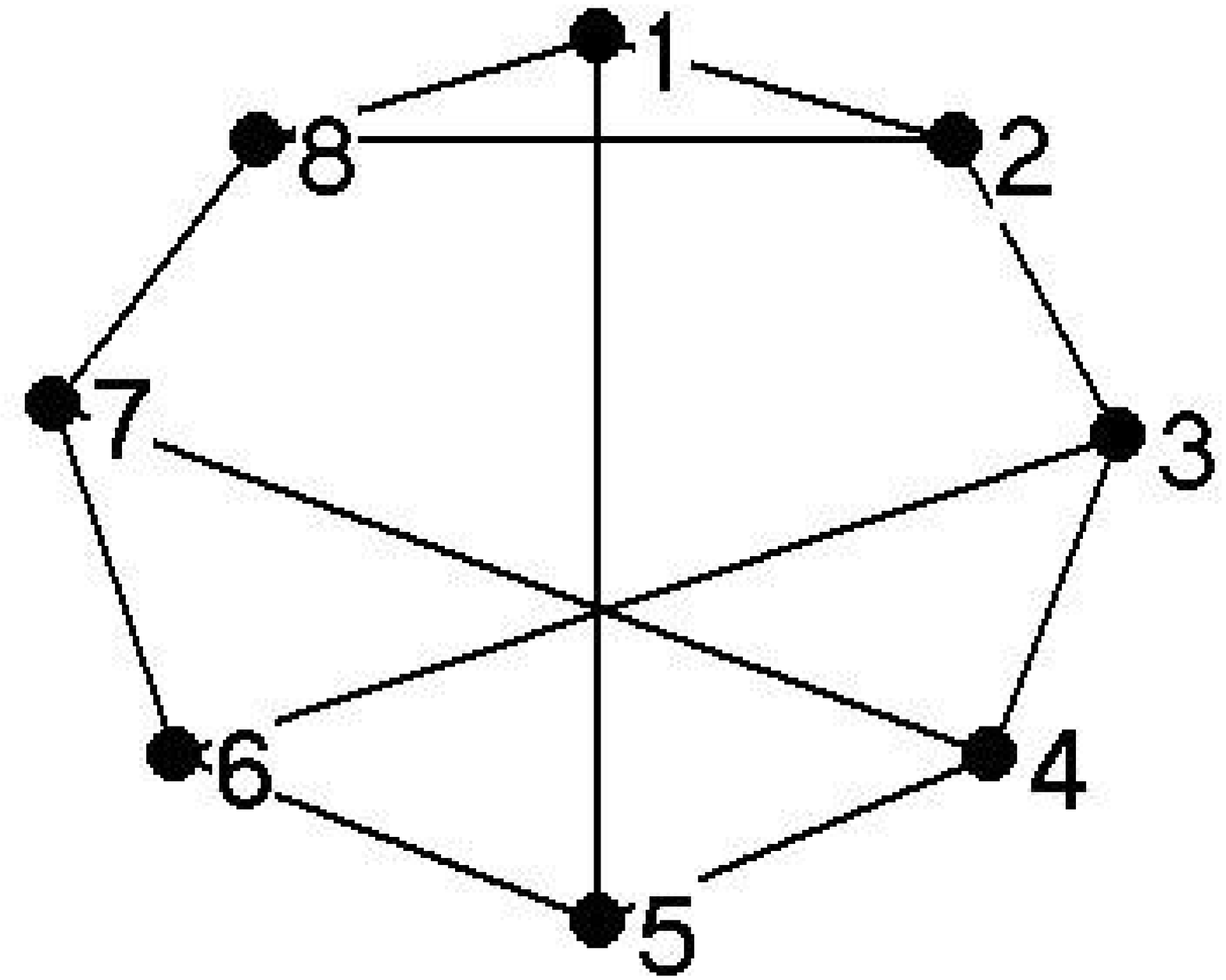}
\includegraphics[angle=90, width=5.5cm]{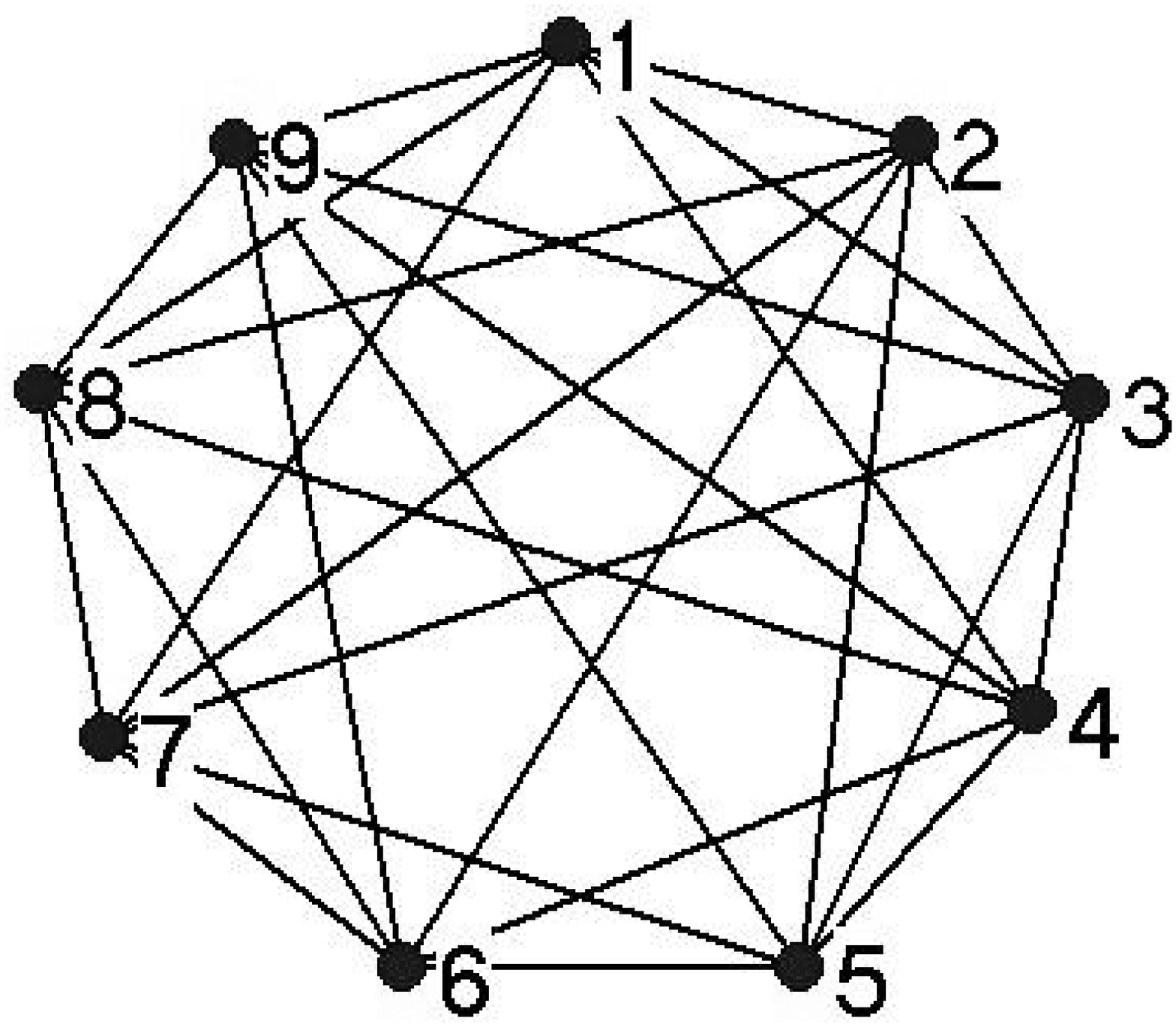}

\end{picture}
\end{center}
\vspace*{-0.8cm}
\end{figure}

In Figure \ref{fig}, we illustrate two regular graphs, with eight
and nine nodes, and degrees equal to 3 and 6, respectively. In the
left hand side graph, nodes $\{1,6,8\}$  are the only ones forming
part of a triangle. Vertices $\{4,6\}$ form part of three squares,
vertices $\{3,5,7\}$ form part of only two squares and the rest do
not form part of any square. The analysis can be obviously extended
to larger subgraphs. However, it is evident that there are three
groups of distinguishable vertices in the graph, $\{1,2,8\}$
,$\{4,6\}$ and $\{3,5,7\}$. These are distinguishable according to
their participation in the different subgraphs, although they cannot
be distinguished by EC. In the right hand side graph, vertices
$\{1,3,5,6,8\}$ take part in 44 of the 100 squares present in the
graph, while vertices $\{2,4,7,9\}$ take part in 45 (all vertices
take part in the same number of smaller subgraphs; e.g., edges,
triangles, connected triples). However, these groups of vertices
cannot be distinguished by any of the centrality measures (DC, CC,
BC and EC).

A method was proposed in \cite{SubgCent} for characterizing nodes in
a network according to the number of closed walks starting and
ending at the node. Closed walks are appropriately weighted such
that their influence on the centrality decreases as the order of the
walk increases:  the contribution of closed walks of length $k$ to
the centrality of the vertex is the number of such closed walks
divided by the factorial of $k$. Each closed walk is associated with
a connected subgraph, which means that this measure "counts" the
times that a node takes part in the different connected subgraphs of
the network, with smaller subgraphs having higher importance.
Consequently,  this measures was called the "subgraph centrality"
(SC) for nodes in a network.

In this paper we propose a generalization of the subgraph
centrality. We propose a general method for characterizing nodes in
the graph according to the number of closed walks starting and
ending at the node. Closed walks are appropriately weighted
according to the topological features that we need to measure.



\section{Functional centrality}

 The following well-known result will be useful in the spectral study of
centralities.

\begin{theorem} {\rm \cite{Biggs}}
Let $v_i$ and $v_j$ be vertices of a  graph $\Gamma$. Let
$\textbf{A}$ be the adjacency matrix of $\Gamma$. Then, the number
of walks of length $k$ in $\Gamma$, from  $v_i$ to $v_j$, is the
entry in position $(i,j)$ of the matrix $\textbf{A}^k$.
\end{theorem}

 Since the adjacency
matrix,  $\textbf{A}$, of $\Gamma$  is a symmetric matrix with real
entries, there exists an orthogonal matrix $U=(u_{ij})$ such that
$\textbf{A}=UDU^T$ where $D=diag(\lambda_1,\lambda_2,...,\lambda_n)$
whose diagonal entries are the eigenvalues of  $\textbf{A}$, and the
columns of $U$ are the corresponding eigenvectors that form an
orthogonal basis of the Euclidean space $\mathbb{R}^n$. Hereafter we
will denote  by $\lambda$ the main eigenvalue of $\textbf{A}$ and
the $j$th column of $U$ will be denoted by $U_j$. It must be
emphasized that, if the graph $\Gamma$ is connected, then the
symmetric and non-negative matrix $\textbf{A}$ is irreducible. As a
consequence, the main eigenvalue of $\textbf{A}$ has a positive
eigenvector of multiplicity one. This fact facilitates the use of
the main eigenvector as a measure of centrality \cite{Bonacich}. It
follows that  walks of length $k$ in $\Gamma$, from $v_i$  to $v_j$,
are
\begin{equation}\label{eq0}
\mu_k(ij)=\left(\textbf{A}^k\right)_{ij}=\sum_{s=1}^n u_{is}
u_{js} \lambda_s^k.\end{equation}

 Moreover, the number of
closed walks of length $k$ starting and ending on vertex  $v_i$ in
$\Gamma$ is given by the local spectral moments $\mu_k(i)$, which
are simply defined as the $i$th diagonal entry of the  $k$th power
of the adjacency matrix, $\textbf{A}$:

\begin{equation} \label{eq1}
 \mu_k(i)=\left(\textbf{A}^k\right)_{ii}=\sum_{s=1}^n\left(u_{is}\right)^2\lambda_s^k.
 \end{equation}

Let $\lambda$ be the spectral radius of $\Gamma$. Let $f$ be a
function, $f:\mathbb{R}\mapsto \mathbb{R}$, whose Taylor series is
$f(x)=\displaystyle\sum_{k=0}^\infty a_kx^k$,
 $\vert x \vert <\lambda_*$, where  $\lambda_*>\lambda$.
 We define the \emph{functional centrality},
$C_{f}(i)$, as
\begin{equation}\label{def}C_{f}(i):=\sum_{j=0}^\infty  a_j\mu_j(i)  . \end{equation}

In this case, the number of closed walks of length $l$ is weighted
by $a_l$. Thus, we can select the function $f$ according to the
topological features that we need to measure. Some interesting cases
will be explained in the next sections.

 For any
$v_i\in V$ we denote by $\ell^{i}(\mathbb{N})$ the space of real
sequences $y=(y_j)_{j=0}^\infty$ such that
$$\sum_{j=0}^\infty  y_j\mu_j(i) <\infty.$$

\begin{theorem}\label{thcentrality}
Let $\Gamma=(V,E)$ be a simple graph of order $n$. Let $\lambda_1\ge
\lambda_2\ge \cdots \ge \lambda_n$ be the eigenvalues of $\Gamma$
and let $U=(u_{ij})$ denotes an orthogonal matrix whose columns are
the corresponding eigenvectors.
 Let $\lambda_*>\lambda_1$ and let $f$ be a real function such that
$f(x)=\displaystyle\sum_{k=0}^\infty a_kx^k$, for $\vert x \vert
<\lambda_*$. Then for each  $v_i\in V$,
$a=(a_0,a_1,...,a_k,...)\in\ell^{i}(\mathbb{N})$ and the functional
centrality is
$$C_{f}(i)=\displaystyle\sum_{j=1}^n(u_{ij})^2 f({\lambda_j}).$$
\end{theorem}

\begin{proof} By the definition of $C_f(i)$ and (\ref{eq1}), we obtain
\begin{equation}\label{reorderin}
C_{f}(i)=\sum_{k=0}^\infty a_k
\left(\sum_{j=1}^n\lambda_j^k(u_{ij})^2\right).
\end{equation}
On the other hand, series (\ref{reorderin}) is obtained by adding
term by term the following convergent series:
\begin{align*}
(u_{i1})^2\displaystyle\sum_{k=0}^\infty a_k\lambda_1^k &=
(u_{i1})^2f(\lambda_1)\\
 (u_{i2})^2\displaystyle\sum_{k=0}^\infty
a_k\lambda_2^k&= (u_{i2})^2f(\lambda_2)\\
&\vdots \\
(u_{in})^2\displaystyle\sum_{k=0}^\infty a_k\lambda_n^k &=
(u_{in})^2f(\lambda_n).
\end{align*}
Thus, series (\ref{reorderin}) converges to
$\displaystyle\sum_{j=1}^n(u_{ij})^2 f({\lambda_j}).$
\end{proof}

\begin{theorem}
Let $\Gamma=(V,E)$ be a nontrivial simple graph of order $n$ and
spectral radius  $\lambda$. Let  $\lambda_*>\lambda$ and let $f$ be
a real function such that $f(x)=\displaystyle\sum_{k=0}^\infty
a_kx^k$, for $\vert x \vert <\lambda_*$, where  $a_l\ge 0$ for all
$l$. Then for each $v_i\in V$,
$$C_{f}(i)\le \frac{1}{n}[f(n-1)+(n-1)f(-1)].$$ The equality holds if and
only if  $\Gamma$ is the complete graph.
\end{theorem}

\begin{proof} Since $\Gamma$ is nontrivial, let $x$ be an edge of $\Gamma$.
Let $\Gamma-x$ be the graph obtained by removing $x$ from $\Gamma$.
Then the number of closed walks of length $k$ in $\Gamma-x$ is equal
to the number of closed walks of length $k$ in $\Gamma$ minus the
number of closed walks of length $k$ in $\Gamma$ containing $x$.
Consequently, as $a_l\ge 0$ for all $l$, for all $v_i\in V$,
$C_f(i)$ in $\Gamma-x$ is lower than or equal to $C_f(i)$ in
$\Gamma$. In closing, the maximum $C_f(i)$ is attained if and only
if $\Gamma$ is the complete graph $\Gamma=K_n$. We now compute
$C_f(i)$ in $K_n$. The eigenvalues of $\Gamma$ are $n-1$ and $-1$
(with multiplicity $1$ and $n-1$). By spectral decomposition of unit
vector $e_i\in \mathbb{R}^n$,
$$e_i=u_{i1}U_1+\sum_{j=2}^nu_{ij}U_j,$$ we obtain
$$1=\parallel e_i\parallel^2=\frac{1}{n}+\sum_{j=2}^n(u_{ij})^2.$$
Therefore, we have
\begin{align*}
C_{f}(i)&=\displaystyle\sum_{j=1}^n(u_{ij})^2 f({\lambda_j})\\
&=\frac{1}{n}f(n-1)+f(-1)\sum_{j=2}^n(u_{ij})^2\\
 &=\frac{1}{n}\left[f(n-1)+(n-1)f(-1)\right] .
\end{align*}
\end{proof}

\subsection{Subgraph Centrality}

 The \emph{subgraph centrality}, \cite{SubgCent}, is defined as
\begin{equation}\label{serie}
C_{S}(i):=\sum_{k=0}^\infty\frac{\mu_k(i)}{k!}.
\end{equation}

In this case, closed walks are appropriately weighted such that
their influence on the centrality decreases as the order of the walk
increases. Each closed walk is associated with a connected subgraph,
which means that this measure "counts" the times that a node takes
part in the different connected subgraphs of the graph, with smaller
subgraphs having higher importance.

Let $\lambda$  be the main eigenvalue of $\textbf{A}$. For any
non-negative integer $k$ and any  $v_i\in V(\Gamma)$ we have
$\mu_k(i)\le \lambda^k$ and thus series (\ref{serie}) - whose terms
are nonnegative - converges:
\begin{equation} \label{serieConverge}
\sum_{k=0}^\infty\frac{\mu_k(i)}{k!}\le
\sum_{k=0}^\infty\frac{\lambda^k}{k!}=e^\lambda .
\end{equation}

Thus,
$a=\left(1,1,\frac{1}{2!},\frac{1}{3!},...,\frac{1}{k!},...\right)\in
\ell^{i}(\mathbb{N})$, $\forall v_i\in V(\Gamma)$. Obviously, in
this case, $f(x)=e^x$ and the subgraph centrality  of $v_i\in V$ is
$C_{S}(i)=\displaystyle\sum_{j=1}^n(u_{ij})^2e^{\lambda_j}.$

The reader is referred to \cite{SubgCent} for a detailed study on
subgraph centrality and its applications to complex networks.

\subsection{Monomial centrality}

Suppose that the problem is to characterize nodes in the graph
according to the number of closed walks of length $k$ containing the
node, i.e., the local spectral moments. In this case we take the
function $f(x)=x^k$ and, in consequence, we call this centrality
\emph{monomial centrality}. Thus, the monomial centrality
$C_{k}(i)$, is defined as
$$C_{k}(i):=\mu_k(i).$$
For instance, if $k=2$, we obtain the \emph{degree centrality}.


\subsection{Functional centrality of radius $k$ }

 Suppose that the problem is to characterize nodes in the graph
according to the number of closed walks of length  lower than or
equal to $k$ containing the node. In this case we take the function
$f$ as a polynomial of degree $k$, $f_k(x)=a_0+a_1x+a_2x^2+\cdots
+a_kx^{k}$, where the coefficients are taking according to the
importance of the participation of the node in the closed walks of a
given length. For instance, if we consider that the contribution of
closed walks of length $1\le l\le k$ to the centrality is inversely
proportional to $l$, then we take $a_l=\frac{1}{l}$. Thus,
$C_{f_k}(i)=\displaystyle\sum_{j=1}^n(u_{ij})^2\left(
1+\lambda_j+\frac{\lambda_j^2}{2}+\cdots +
 \frac{\lambda_j^{k}}{k}\right)
 $.

\subsection{Odd and even centralities}

It is well-known that there are graphs that do not have odd
closed-walks. i.e., the bipartite graphs. On the other hand, it
would be of some interest to characterize nodes in the graphs
according to the number  of closed walks of odd (even) length
containing the node. In such a case we can take an odd (even)
function. For instance, the \emph{odd subgraph centrality} is
defined as
$$
C_{S_{odd}}(i):=\frac{\mu_1(i)}{1!}+\frac{\mu_3(i)}{3!}+\frac{\mu_5(i)}{5!}+\cdots
$$
Hence, in this case, $f(x)=sinh(x)$ and
$C_{S_{odd}}(i)=\displaystyle\sum_{j=1}^n(u_{ij})^2sinh({\lambda_j}).$
Analogously, the \emph{even subgraph centrality} is defined as
$$
C_{S_{even}}(i):=\frac{\mu_0(i)}{0!}+\frac{\mu_2(i)}{2!}+\frac{\mu_4(i)}{4!}+\cdots
$$
Thus,  $f(x)=cosh(x)$ and
$C_{S_{even}}(i)=\displaystyle\sum_{j=1}^n(u_{ij})^2cosh({\lambda_j}).$

As we will show in Section \ref{bipartivity}, the above centralities
 are of singular importance in the study of bipartivity.

\subsection{Example}


In order to show the differences in the orders imposed by different
functional centralities, we have selected the left hand side graph
of Figure \ref{fig} and the following centrality measures.
\begin{itemize}
 \item Monomial centrality, $C_3$: $f(x)=x^3$.

 \item Functional centrality of radius 3, $C_{f_3}$: $f(x)=1+x+\frac{x^2}{2}+\frac{x^3}{3}$.

 \item Functional centrality of radius 4, $C_{f_4}$: $f(x) = 1+x+\frac{x^2}{2}+\frac{x^3}{3}+\frac{x^4}{4}$.

 \item Odd centrality, $C_{S_{odd}}$: $f(x)=\sinh(x)$.

  \item Subgraph centrality, $C_s$: $f(x)=e^x$.

 \end{itemize}

{\small
\begin{center}
\begin{tabular}{|c|c|c|c|c|c|c|c|c|c|}
\hline
& {\bf 1 }& {\bf 2 } & {\bf 3} & {\bf 4} & {\bf 5}& {\bf 6} & {\bf 7 }& {\bf 8}  \\
\hline
$C_3$  & 2 & 2 &  0 &0 & 0& 0 & 0 & 2\\
\hline
$C_{f_3}$ &3.16 & 3.16  & 2.5  &   2.5  &  2.5 &   2.5& 2.5 & 3.16 \\
\hline
$C_{f_4}$ &6.92& 6.92 & 7.25 &  7.75 &  7.25 & 7.75 & 7.25 & 6.92  \\
\hline
$C_{S_{odd}}$ & 0.608&0.608 & 0.117 &0.075  &0.117 & 0.075 & 0.117  & 0.608  \\
\hline
$C_s$ & 3.9 & 3.9 & 3.63  & 3.7 &3.63 & 3.7 & 3.63 & 3.9\\
\hline
\end{tabular}
\end{center}
}

\section{Functional centralization}\label{bipartivity}

By (\ref{eq0}) and (\ref{eq1}) we have that the number $w_k$ of
walks of length $k$ in $\Gamma$ is given by
\begin{equation}w_k=\sum_{i,j}\mu_k(ij)=\sum_{s=1}^n\left(\sum_{i=1}^n
u_{is}\right)^2\lambda_s^k,\end{equation} and the number $\theta_k$
of closed-walks of length  $k$ in  $\Gamma$ is given by the trace of
$\textbf{A}^k$:
\begin{equation}
\theta_k=\displaystyle \sum_{i=1}^n
C_{k}(i)=\displaystyle\sum_{s=1}^n\lambda_s^k .
\end{equation}

Let $\lambda$ be the spectral radius of $\Gamma$. Let $f$ be a
function, $f:\mathbb{R}\mapsto \mathbb{R}$, whose Taylor series is
$f(x)=\displaystyle\sum_{k=0}^\infty a_kx^k$,
 $\vert x \vert <\lambda_*$, where  $\lambda_*>\lambda$.
 We define the \emph{functional centralization},
$C_{f}(\Gamma)$, as
\begin{equation}C_{f}(\Gamma):=\sum_{j=0}^\infty  a_j\theta_j  . \end{equation}

Hence, an analytical expression for $C_{f}(\Gamma)$, which depends
only on the eigenvalues and order of the graph, can be obtained by
using a procedure analogous to that described in the proof of
Theorem \ref{thcentrality}.

\begin{theorem}\label{thGcentrality}
Let $\Gamma=(V,E)$ be a simple graph of order $n$ and let
$\lambda_1\ge \lambda_2\ge \cdots \ge \lambda_n$ be the eigenvalues
of $\Gamma$. Let $\lambda_*>\lambda_1$ and let $f$ be a real
function such that $f(x)=\displaystyle\sum_{k=0}^\infty a_kx^k$, for
$\vert x \vert <\lambda_*$. Then
 the
functional centralization is
$$C_{f}(\Gamma)=\displaystyle\sum_{j=1}^n
f({\lambda_j})=\displaystyle\sum_{i=1}^n C_{f}(i).$$
\end{theorem}

For instance, taking $f(x)=x^k$ we obtain the $k$-\emph{spectral
moment} of $\Gamma$,
$$C_{k}(\Gamma):=\theta_k=\displaystyle\sum_{s=1}^n\lambda_s^k .$$
This centralization
 only measures  the number of closed-walks of length $k$ in $\Gamma$.  For instance,
 if $k=2$, then $C_{2}( G)=2m$, where $m$ denotes the size of  $G$.

 Taking $f(x)=e^x$ we obtain the \emph{subgraph centralization}:
\begin{equation}
C_{S}(\Gamma):=\sum_{k=0}^\infty \frac{\theta_k}{k!}.
\end{equation}
That is,
\begin{equation}
 C_{S}(\Gamma)=\sum_{i=1}^n e^{\lambda_i}.
\end{equation}

Taking $f(x)=\sinh(x)$ in Theorem \ref{thGcentrality} we obtain the
\emph{odd subgraph centralization} of $\Gamma$
\begin{equation}
C_{S_{odd}}(\Gamma)=\displaystyle\sum_{i=1}^nsinh({\lambda_i}).
\end{equation}
Analogously, the \emph{even subgraph centralization} is
\begin{equation}
C_{S_{even}}(\Gamma)=\displaystyle\sum_{i=1}^n cosh({\lambda_i}).
\end{equation}

Both subgraph centralization and even centralization, were used to
obtain measures of bipartivity in graphs \cite{bipart}: the
proportion of even closed walks to the total number of closed walks
is a measure of the graph bipartivity. Thus, we measures the
bipartivity of a graph $\Gamma$ by
$$\beta(\Gamma)=\frac{\displaystyle\sum_{i=1}^n\cosh(\lambda_i)}{\displaystyle\sum_{i=1}^ne^{\lambda_i}}.$$
It is evident that  $\beta(\Gamma)\le 1$ and $\beta(\Gamma)= 1$ if
and only if $\Gamma$ is bipartite. Moreover, as
$e^{\lambda_i}=\sinh(\lambda_i)+\cosh(\lambda_i)$ and
$\sinh(\lambda_i)< \cosh(\lambda_i)$, $\forall i$, we have
$\frac{1}{2}<\beta(\Gamma)\le 1$.  The lower bound is reached for
the least possible bipartite graph with $n$ nodes, which is the
complete graph $K_n$. As the eigenvalues of $K_n$ are $n-1$ and $-1$
(with multiplicity $n-1$), then $\beta(\Gamma)\rightarrow
\frac{1}{2}$ when $n\rightarrow \infty$ in $K_n$.

The contribution of vertex $v_i$ to graph bipartivity, $\beta(i)$,
can be obtained by using the appropriate functional centrality of
$v_i$. That is,
$$\beta(i)=\frac{\displaystyle\sum_{j=1}^n (u_{ij})^2\cosh(\lambda_j)}
{\displaystyle\sum_{j=1}^n (u_{ij})^2e^{\lambda_j}}.$$

The reader is referred to \cite{bipart} for a more detailed study on
bipartivity and its applications. Moreover, an application of our
bipartivity measures to the study of fullerene graphs can be found
in \cite{fullerene}.

\section{Functional centrality in hypergraphs}

Let ${\cal  H}$  be a simple hypergraph of order $n$. The adjacency
matrix, $\textbf{A}({\cal H})$, of the hypergraph ${\cal H}=(V,E)$
is a square symmetric matrix whose entries $a_{ij}$ are the number
of hyper-edges that contain both nodes $v_i$ and $v_j$: the diagonal
entries of  $\textbf{A}({\cal H})$ are zero. This can be obtained
from the incidence matrix of ${\cal H}$ as follows:
$$ \textbf{A}({\cal
H})=\textbf{E} \textbf{E}^T-\textbf{D}$$ where $\textbf{E}^T$  is
the transpose of the incidence matrix and $\textbf{D}$ is the
diagonal matrix whose diagonal entries are the degrees of the
vertices. More formally, $\textbf{A}({\cal H})$  is a $n\times n$
matrix with diagonal entries $a_{ii}=0$, for $v_i\in V$, and
off-diagonal entries
$$a_{ij}=|\{E_k\in
E({\cal  H}):\{v_i,v_j\}\subset E_k\}|, \mbox{\rm { } for { }}
v_i,v_j\in V({\cal  H}),\quad i\neq j.$$ Since $\textbf{A}({\cal
H})$ is symmetric, and its entries are non-negative integers, it may
be viewed as the adjacency matrix of a multigraph  $G'$, i.e., a
graph having multiple links between nodes, called the
\emph{associated graph} of ${\cal  H}=(V,E)$.

It must be emphasized that, if the hypergraph ${\cal  H}$ is
connected, then the symmetric and non-negative matrix
$\textbf{A}({\cal  H})$ is irreducible. As a consequence, the main
eigenvalue of $\textbf{A}({\cal  H})$ has a positive eigenvector of
multiplicity one. This fact facilitates the extension, to the case
of hypergraphs, of the use of the main eigenvector as a measure of
centrality.

The following result   will be useful in extending the definition of
functional centrality to hypergraphs.

\begin{theorem} {\rm \cite{Walk-Reg-Hyp}}
Let $v_i$ and $v_j$ be vertices of a  hypergraph ${\cal  H}$. Let
$\textbf{A}({\cal  H})$ be the adjacency matrix of ${\cal  H}$.
Then, the number of walks of length $k$ in ${\cal  H}$, from  $v_i$
to $v_j$, is the entry in position $(i,j)$ of the matrix
$\left(\textbf{A}({\cal H})\right)^k$.
\end{theorem}

So, the spectral study of centralities is completely analogous to
the previous one for graphs.

\begin{theorem}\label{thcentrality1}
Let ${\cal H}$ be a simple hypergraph of order $n$. Let
$\lambda_1\ge \lambda_2\ge \cdots \ge \lambda_n$ be the eigenvalues
of $G_{\cal H}$ and let $\textbf{U}=(\textbf{u}_{ij})$ denote an
orthogonal matrix whose columns are the corresponding eigenvectors.
Let $\lambda_*>\lambda_1$ and let $f$ be a real function such that
$f(x)=\displaystyle\sum_{k=0}^\infty a_kx^k$, for $\vert x \vert
<\lambda_*$. Then for each  $v_i\in V({\cal H})$,
$a=(a_0,a_1,...,a_k,...)\in\ell^{i}(\mathbb{N})$ and the functional
centrality  is
$$C_{f}(i)=\displaystyle\sum_{j=1}^n(\textbf{u}_{ij})^2
f({\lambda_j}).$$
\end{theorem}

\begin{proof}
The proof is basically as in Theorem \ref{thcentrality}.
\end{proof}

The particular case of subhypergraph centrality, i.e., $f(x)=e^x$,
 was analyzed in \cite{hypergraphs}.

\end{document}